\documentclass{article} 
\usepackage{amssymb}
\usepackage{amsmath}
\usepackage{euscript}
\begin{document} 
\input epsf.sty

\title{Algebraic curves $P(x)-Q(y)=0$ and 
functional equations}
 
\author{F. Pakovich}

\date{} 

\maketitle

\def\be{\begin{equation}}
\def\ee{\end{equation}}
\def\bs{$\square$ \vskip 0.2cm}
\def\d{{\rm d}} 
\def\D{{\rm D}} 
\def\I{{\rm I}} 
\def\C{{\mathbb C}} 
\def\N{{\mathbb N}} 
\def\P{{\mathbb P}}
\def\Z{{\mathbb Z}}
\def\R{{\mathbb R}} 
\def\ord{{\rm ord}}
\def\ssigma{\omega}
\def\f{\EuScript}
\def\e{\eqref}
\def\phi{{\varphi}}
\def\v{{\varepsilon}} 
\def\deg{{\rm deg\,}} 
\def\Det{{\rm Det}}
\def\dim{{\rm dim\,}} 
\def\Ker{{\rm Ker\,}} 
\def\Gal{{\rm Gal\,}}
\def\card{{\rm card}}
\def\St{{\rm St\,}} 
\def\exp{{\rm exp\,}} 
\def\cos{{\rm cos\,}} 
\def\diag{{\rm diag\,}} 
\def\GCD{{\rm GCD }}
\def\LCM{{\rm LCM }}
\def\mod{{\rm mod\ }}
\def\Res{{\rm Res\ }}

\def\bp{\begin{proposition}}
\def\ep{\end{proposition}}
\newtheorem{zzz}{Theorem}
\newtheorem{yyy}{Corollary}
\def\bt{\begin{theorem}}
\def\et{\end{theorem}}
\def\be{\begin{equation}}
\def\bee{\begin{equation*}}
\def\la{\label}
\def\l{\lambda}
\def\m{\mu}
\def\ee{\end{equation}}
\def\eee{\end{equation*}}
\def\bl{\begin{lemma}}
\def\el{\end{lemma}}
\def\bc{\begin{corollary}}
\def\ec{\end{corollary}}
\def\proof{\noindent{\it Proof. }}
\def\note{\noindent{\bf Note. }}
\def\bd{\begin{definition}}
\def\ed{\end{definition}}
\def\e{\eqref}
\def\qed{$\ \ \Box$ \vskip 0.2cm}
\newtheorem{theorem}{Theorem}[section]
\newtheorem{lemma}{Lemma}[section]
\newtheorem{definition}{Definition}[section]
\newtheorem{corollary}{Corollary}[section]
\newtheorem{proposition}{Proposition}[section]

\begin{abstract}
In this paper we give several conditions implying the irreducibility of the algebraic curve
$P(x)-Q(y)=0,$ 
where $P,Q$ are rational functions. We also
apply the results obtained to the functional equations $P(f)=Q(g)$ \linebreak and $P(f)=cP(g),$ where $c\in \C.$  
For example, we show that 
for a generic pair of rational functions $P,Q$ 
the first equation has no non-constant solutions 
$f,g$ meromorphic on $\C$
whenever $(\deg P-1)(\deg Q-1)\geq 2.$ 
\end{abstract}

\section{Introduction}  In the paper \cite{hy} K. H. Ha and C. C. Yang proved that if $P,Q$ is a pair of polynomials
such that $P$ and $Q$ have no common finite critical values 
and 
$n=\deg P$ and $m= \deg Q$
satisfy some constraints then the functional equation \be \la{ur} P(f)=Q(g)\ee has no non-constant solutions $f,g$ meromorphic on $\C.$
This result yields in particular 
that for given $n,m$ satisfying above constraints
there exists a proper algebraic subset 
$\Sigma\subset \C^{n+m+2}$ such that for any pair of polynomials
$$P(z)=a_nz^n+a_{n-1}z^{n-1}+...+a_1z+a_0, \ \ \ Q(z)=b_mz^m+a_{m-1}z^{m-1}+...+b_1z+b_0$$ 
with $(a_n,...,a_0,b_m, ... , b_0)\notin \Sigma$ equation \eqref{ur} 
has no non-constant solutions $f,g$ meromorphic on $\C.$
Some further results concerning equation \eqref{ur} were obtained in the papers \cite{es0}, \cite{ly}, \cite{ly2}, 
\cite{pak2}.

The approach of \cite{hy} is based on the Picard theorem which states that an algebraic curve $q(x,y)=0$ 
of genus $\geq 2$ can not be parametrized by non-constant functions $f,g$ 
meromorphic on $\C$. The Picard theorem implies that for given polynomials $P,Q$ equation \eqref{ur} has non-constant meromorphic solutions $f,g$ if and only if the algebraic curve \be \la{cur0} P(x)-Q(y)=0 \ee has an irreducible component of genus $\leq 1.$ Indeed, any non-constant solution $f,g$ of \eqref{ur} parametrizes an irreducible
component of \eqref{cur0} and the genus of this component equals 0 or 1 by the Picard theorem.  
On the other hand, any irreducible component of genus
0 or 1 of curve \eqref{cur0} may be parametrized correspondingly by  
rational or elliptic functions $f,g$. Clearly, these functions satisfy \eqref{ur} and hence \eqref{ur} has meromorphic solutions.

A question closely related to equation \eqref{ur} is the problem of description of 
so called ``strong uniqueness polynomials'' for meromorphic functions that is of polynomials $P$ such that 
the equality 
\be \la{c} P(f)=cP(g)\ee 
for $c\in \C$ and non-constant functions $f,g$ meromorphic on $\C$ 
implies that $c=1$ and $f\equiv g.$ This problem arose in  
connection with the problem of description of ``uniqueness range sets" for meromorphic functions and 
was studied in the recent papers 
\cite{hy,ly,ly2,pak2,a0,a,es,cher,fu1,fu2, fu3, sh, yh}.
Clearly, the Picard theorem is 
applicable to this problem too. Namely, it follows from the Picard theorem that $P$ is a strong uniqueness polynomial for meromorphic functions
if and only if for any $c\neq 1$ the curve $P(x)-cP(y)=0$ has no irreducible components of genus $\leq 1$, and a unique such a component of the curve $P(x)-P(y)=0$ is $x-y=0$ (the last condition is obviously 
equivalent to the condition that 
the curve 
\be \la{kot}
\frac{P(x)-P(y)}{x-y}=0
\ee
has no irreducible components of genus $\leq 1$).

Although the Picard theorem reduces the question about the existence of meromorphic solutions of equation \eqref{ur}
to an essentially algebraic question about curve \eqref{cur0} 
most of the papers concerning
equation \eqref{ur} or strong uniqueness polynomials for meromorphic functions 
use the Nevanlinna value distribution theory and other analytic methods. Actually, 
the algebraic methods seem to be underestimated and one of the goals of this
paper is to show that these methods are 
not less fruitful and sometimes lead to more precise results than the analytic ones. 

In this paper we consider equations \eqref{ur}, \eqref{c} for arbitrary {\it rational} $P$ and $Q$ and show that 
for ``generic'' $P,Q$ they have only ``trivial'' meromorphic solutions whenever the degrees of $P$ and $Q$ satisfy some 
mild restrictions. 
It is easy to see that the Picard theorem is still applicable to equations \eqref{ur} and \eqref{c} with rational $P,Q$ 
if instead of curves \eqref{cur0} and \eqref{kot} 
to consider correspondingly the curves 
\be \la{cur} h_{P,Q}(x,y):\ P_1(x)Q_2(y)-P_2(x)Q_1(y)=0,\ee 
and 
\be \la{cur2} h_P(x,y):\ \frac{P_1(x)P_2(y)-P_2(x)P_1(y)}{x-y}=0\ee 
where 
$P_1$, $P_2$ and $Q_1,Q_2$ are pairs polynomials without common roots such that 
$P=P_1/P_2,$ $Q=Q_1/Q_2.$
An explicit description of pairs of rational functions $P,Q$ for which the curve $h_{P,Q}(x,y)$
is irreducible is known only in the case when $P,Q$ are indecomposable polynomials (see \cite{couv}).
On the other hand, in order to analyze equations \eqref{ur} and \eqref{c} for generic
rational functions $P,Q$ it is necessary to have available conditions implying the irreducibility of 
curves \eqref{cur} 
and \eqref{cur2}
for wide classes of $P,Q$.
In this paper, using the description of irreducible components of \eqref{cur}
given in \cite{pak},
we provide several such conditions and apply the results obtained to 
equations \eqref{ur} and \eqref{c}. Recall that a point $s\in \C\P^1$ is called a critical value of a rational 
function $F$ if the set $F^{-1}\{s\}$ contains 
less than $\deg F$ points, and $s$ is called a simple critical value if $F^{-1}\{s\}$ contains 
exactly $\deg F-1$ points. We will denote the set of all critical values of $F$ by $\f C(F).$ 

Our main result concerning curve \eqref{cur} is a complete analysis
of its irreducibility in the case when $\f C(P)\cap \f C(Q)$ contains ``few'' elements. 
Namely, 
we show 
that curve \eqref{cur} is irreducible whenever $\f C(P)\cap \f C(Q)$ is empty or contains one point and give an explicit condition for its irreducibility in the case when $\f C(P)\cap \f C(Q)$ contains two points. 
Besides, we show 
that curve \eqref{cur2} is irreducible if $P$ is indecomposable 
and has at least one simple critical value, or if all critical values of $P$ are simple.

As an application of our results about curves \eqref{cur} and \eqref{cur2} we obtain several results  
concerning equations \eqref{ur} and \eqref{c}.
In particular, we prove analogues of the results of \cite{hy} for rational $P,Q$. Our main result concerning equation \eqref{ur} is the following theorem. 

\bt \la{1} Let $P,Q$ be a pair of rational functions
such that $\f C(P)\cap \f C(Q)=\emptyset$.
Then functional equation \eqref{ur} has  
non-constant solutions $f,g$ meromorphic on $\C$ if and only if   
$n=\deg P$ and $m=\deg Q$ satisfy the inequality
$(n-1)(m-1)< 2$. 
\et 

From Theorem \ref{1} we deduce the following result.\footnote{Professor C. C. Yang kindly informed us that 
a similar result is obtained by a different method also in the forthcoming paper \cite{forth}.}

\bt \la{ca1} Let $n,m$ by any integer non-negative numbers such that the inequality $(m-1)(n-1)\geq 2$ holds.
Then there exists a proper algebraic subset 
$\Sigma\subset \C\P^{2n+1}\times\C\P^{2m+1}$ such that for any pair of rational functions
$$P(z)=\frac{a_nz^n+a_{n-1}z^{n-1}+...+a_0}{b_nz^n+b_{n-1}z^{n-1}+...+b_0}, \ \ \ Q(z)=\frac{c_mz^m+c_{m-1}z^{m-1}+...+c_0}{d_mz^m+d_{m-1}z^{m-1}+...+d_0}$$ 
with $(a_n,...,a_0,b_n, ... , b_0,
c_m,...,c_0,d_m, ... , d_0)\notin \Sigma$ equation \eqref{ur} has no  
non-cons\-tant solutions $f,g$ meromorphic on $\C$. 
\et

Furthermore, we prove an analogue of Theorem \ref{1} for the functional equation \be \la{ur2} P(f)=P(g),\ee where $P$ is a rational function,
generalizing the previous result of the paper \cite{sh} concerning the case when $P$ is a polynomial.

\bt \la{2} Let $P$ be a
rational function of degree $n$ which has only 
simple critical values. Then functional equation \eqref{ur2} has  
non-constant solutions $f,g$ such that $f\not\equiv g$ and $f,g$ are meromorphic on $\C$ if and only if 
$n< 4.$ 
\et

Finally, from Theorems \ref{1} and \ref{2} we deduce the following theorem. 

\bt \la{ca2} For any $n\geq 4$ there exists a proper algebraic subset 
$\Sigma\subset \C\P^{2n+1}$ such that for any rational function
$$P(z)=\frac{a_nz^n+a_{n-1}z^{n-1}+...+a_0}{b_nz^n+b_{n-1}z^{n-1}+...+b_0}$$ 
with $(a_n,...,a_0,b_n, ... , b_0)\notin \Sigma$ equality \eqref{c}, where $f,g$ are non-constant 
functions meromorphic on $\C$, implies that $c=1$ and $f\equiv g$. 
\et

The paper is organized as follows. In the second section we recall a construction from \cite{pak} which permits to 
describe irreducible components of \eqref{cur} and \eqref{cur2} and 
to calculate their genuses. 
In the third section we give several conditions implying the irreducibility of curves \eqref{cur} and \eqref{cur2}.
Finally, in the fourth section we prove our results concerning equations \eqref{ur} and \eqref{c}.

\section{Components of $h_{P,Q}(x,y)$ and $h_{P}(x,y)$} 
In this section we recall a construction from \cite{pak} which permits to describe irreducible 
components of the curves $h_{P,Q}(x,y)$ and $h_{P}(x,y)$.

For rational functions $P$ and $Q$  
denote by $S=\{z_1, z_2, \dots , z_r\}$ the union of $\f C(P)$ and $\f C(Q)$.
Fix a point 
$z_0$ from $\C\P^1\setminus S$ and small loops $\gamma_i$ around $z_i,$ $1\leq i \leq r,$ such that 
$\gamma_1\gamma_2...\gamma_r=1$ in $\pi_1(\C\P^1\setminus S,z_0)$. Set $n=\deg P,$ $m=\deg Q$.
For $i,$ $1\leq i \leq r,$ denote by  
$\alpha_i\in S_n$ (resp. $\beta_i\in S_m$) a permutation of points of 
$P^{-1}\{z_0\}$ (resp. of $Q^{-1}\{z_0\}$) induced by the lifting of $\gamma_i$ by $P$ (resp. $Q$).
Clearly, the permutations $\alpha_i$ (resp. $\beta_i$), $1\leq i \leq r,$ generate the monodromy group of 
$P$ (resp. of $Q$) and \be \la{per} \alpha_1\alpha_2...\alpha_r=1,\ \ \ \ \beta_1\beta_2...\beta_r=1.\ee
Notice that since $S=\f C(P)\cup \f C(Q)$ some of permutations $\alpha_i,$ $\beta_i,$ $1\leq i \leq r,$
may be identical permutations.

Define now permutations $\delta_1, \delta_2, \dots, \delta_{r}\in S_{nm}$ as follows:
consider the set of $mn$ elements $c_{j_1,j_2},$ $1\leq j_1 \leq n,$
$1\leq j_2\leq m,$ and set $(c_{j_1,j_2})^{\delta_i}=c_{j_1^{\prime},j_2^{\prime}},$ where $$j_1^{\prime}
=j_1^{\alpha_i},\ \ \ \ j_2^{\prime}=j_2^{\beta_i}, \ \ \ \ 1\leq i \leq r.$$ 
It is convenient to consider $c_{j_1,j_2},$ $1\leq j_1 \leq n,$
$1\leq j_2\leq m,$ as elements of a $n\times m$ matrix $M$. Then the action of the 
permutation $\delta_i,$ $1\leq i \leq r,$ reduces 
to the permutation of rows of $M$ in accordance with the permutation $\alpha_i$ and 
the permutation of columns of $M$ in accordance with the permutation $\beta_i$.

In general, 
the permutation group $\Gamma(P,Q)$ generated by $\delta_i,$ $1\leq i \leq r,$
is not transitive on 
the set $c_{j_1,j_2},$ $1\leq j_1 \leq n,$
$1\leq j_2\leq m$. Denote by $o(P,Q)$ the number of transitivity sets of 
the group $\Gamma(P,Q)$ and let $\delta_{i}(j),$ 
$1\leq i \leq r,$ $1\leq j \leq o(P,Q),$ be the  
permutation induced by the permutation $\delta_i,$ $1\leq i \leq r,$
on the transitivity set $U_j,$ $1\leq j \leq o(P,Q).$ We will denote the permutation group generated by the permutations $\delta_{i}(j),$ $1\leq i \leq r,$ for some fixed $j,$ $1\leq j \leq o(P,Q),$
by $G_j.$ 

By construction, the group $G_j,$ $1\leq j \leq o(P,Q),$
is a transitive permutation group on $U_j.$ 
Furthermore, it follows from \eqref{per} that $\delta_1\delta_2...\delta_r=1$ and hence for any $j,$ 
$1\leq j \leq o(P,Q),$ the equality 
$$\delta_{1}(j)\delta_{2}(j)\dots \delta_{r}(j)=1$$
holds. By the Riemann existence theorem (see e.g. \cite{mir}, Corollary 4.10)
this implies that there exist compact Riemann surfaces $R_j$ and holomorphic functions 
$h_j:\, R_j\rightarrow \C\P^1,$ $1\leq j \leq o(P,Q),$ non ramified outside of $S,$  
such that the permutation $\delta_{i}(j),$ 
$1\leq i \leq r,$ $1\leq j \leq o(P,Q),$
is induced by the lifting of $\gamma_i$ by $h_j.$

Moreover, it follows from the construction of the group $\Gamma(P,Q)$ that for each $j,$ $1\leq j \leq o(P,Q),$
the intersections of the transitivity set $U_j$ with the rows of $M$
form an imprimitivity system $\Omega_P(j)$ for the group $G_j$ such that the permutations
of blocks of $\Omega_P(j)$ induced by 
$\delta_i(j),$ $1\leq i \leq r,$ coincide with $\alpha_i.$ 
Similarly, the intersections of $U_j$ with the columns of $M$
form an imprimitivity system $\Omega_Q(j)$ such that the permutations
of blocks of $\Omega_Q(j)$ induced by 
$\delta_i(j),$ $1\leq i \leq r,$ coincide with $\beta_i.$ 
This implies that  
there exist holomorphic
functions $u_j:\,  R_j\rightarrow \C\P^1$ and $v_j:\,  R_j\rightarrow \C\P^1$
such that 
\be \la{e4} h_j=P\circ u_j=Q\circ v_j,\ee
where the symbol $\circ$ denotes the superposition of functions, $f_1\circ f_2=f_1(f_2)$.

Finally, notice that for any choice of points 
$a\in P^{-1}\{z_0\}$ and $b\in Q^{-1}\{z_0\}$ there exist uniquely defined $j,$ $1\leq j \leq o(P,Q),$ and $c\in h_j^{-1}\{z_0\}$
such that 
\be \la{s}  u_j(c)=a, \ \ \ v_j(c)=b.\ee
Indeed, it is easy to see that if $l,$ $1\leq l \leq n,$ is the index which corresponds to the point $a$  
under the identification of the set $P^{-1}\{z_0\}$ with the set of rows of $M,$ and
$k,$ $1\leq k \leq m,$ is the index which corresponds to the point $b$ 
under the identification of the set $Q^{-1}\{z_0\}$
with the set of columns of $M$, then the needed index $j$ is defined by the condition that
the transitivity set $U_{j}$ contains the element $c_{l,k}$, and 
the needed point $c$ is defined by the condition that 
$c$ corresponds to $c_{l,k}$ under the identification of the set $h_j^{-1}\{z_0\}$ with the set of elements of $U_j.$

\bp[(\cite{pak})]\la{iop}  The Riemann surfaces $R_j,$ $1\leq j \leq o(P,Q),$
are in a one-to-one correspondence with irreducible components of the curve $h_{P,Q}(x,y).$
Furthermore, each $R_j$ is a desingularization of the corresponding component.
In particular, the curve $h_{P,Q}(x,y)$ is irreducible if and only if the group $\Gamma(P,Q)$ is transitive.
\ep

\begin{proof} For $j,$ $1\leq j \leq o(P,Q),$ 
denote by $S_j$ the union of poles of $u_j$ and $v_j$ and 
define the mapping $t_j:\, R_j\setminus S_j \rightarrow \C^2$  
by the formula $$z\rightarrow (u_j(z),v_j(z)).$$
It follows from formula \eqref{e4} that for each $j,$ $1\leq j \leq o(P,Q),$ the mapping 
$t_j$ 
maps $R_j$ to an irreducible component of the curve $h_{P,Q}(x,y)$. Furthermore, for any 
point $(a,b)$ on $h_{P,Q}(x,y)$, such that $z_0=P(a)=Q(b)$ is not contained in $S$, 
there exist uniquely defined $j,$ $1\leq j \leq o(P,Q),$ and $c\in h_j^{-1}\{z_0\}$
satisfying \eqref{s}. This implies that the Riemann surfaces $R_j,$ $1\leq j \leq o(P,Q),$
are in a one-to-one correspondence with irreducible components of $h_{P,Q}(x,y)$ 
and that each mapping $t_j,$ $1\leq j \leq o(P,Q),$
is generically injective. Since an injective mapping of Riemann surfaces is an isomorphism onto an open subset 
we conclude that 
each $R_j$ is a desingularization of the corresponding component of $h_{P,Q}(x,y)$. \end{proof}

For $i,$ $1\leq i \leq r,$ denote by
$$
\l_i=(p_{i,1},p_{i,2}, ... , p_{i,u_i})$$ 
the collection of lengths of disjoint cycles in the permutation $\alpha_i$,
by 
$$
\m_i=(q_{i,1},q_{i,2}, ... , q_{i,v_i})
$$ 
the collection of lengths of disjoint cycles in the permutation $\beta_i$
and by $e_i(j),$ $1\leq i \leq r,$ $1\leq j \leq o(P,Q),$
the number of disjoint cycles in the permutation 
$\delta_i(j)$. 
The Riemann-Hurwitz formula implies that for the genus $g_j,$ $1\leq j \leq o(P,Q),$
of the component of $h_{P,Q}(x,y)$
corresponding to $R_j$ we have:
$$2-2g(R_j)=
\sum_{i=1}^{r}e_i(j)-
\card\{U_j\}(r-2).$$
On the other hand it follows easily from the definition that 
the permutation $\delta_i,$ $1\leq i \leq r,$  
contains  
$$\sum_{j_1=1}^{u_{i}}\sum_{j_2=1}^{v_{i}} \GCD(p_{i,j_1}q_{i,j_2})$$ 
disjointed cycles. In particular, in the case when the curve $h_{P,Q}(x,y)$ is irreducible 
we obtain the following formula for its genus established earlier in \cite{f3}. 

\bc \la{p3} If the curve $h_{P,Q}(x,y)$ is irreducible then for its genus $g$ 
the following formula holds:
\be \la{rh0} 2-2g=
\sum_{i=1}^{r}
\sum_{j_1=1}^{u_{i}} \sum_{j_2=1}^{v_{i}} \GCD(p_{i,j_1}q_{i,j_2})-(r-2)nm. \ \ \ \Box
\ee
\ec

Similarly, we obtain the following corollary concerning the curve $h_{P}(x,y).$

\bc \la{doub} The curve $h_{P}(x,y)$ is irreducible if and only if the monodromy group $G(P)$ of $P$ is 
doubly transitive. Furthermore, if $h_{P}(x,y)$ is irreducible then 
for its genus $g$ 
the following formula holds:
\be \la{rh2}4-2g=
\sum_{i=1}^{r}
\sum_{j_1=1}^{u_{i}} \sum_{j_2=1}^{v_{i}} \GCD(p_{i,j_1}p_{i,j_2})-(r-2)n^2. 
\ee
\ec
\begin{proof} Indeed, it follows from Proposition \ref{iop} that $h_{P}(x,y)=0$ is irreducible if and only if
the group $\Gamma(P,P)$ has two transitivity sets on $M$: the diagonal 
$$\Delta:\ \{c_{j,j} \ \vert \  1\leq j \leq n\}$$ (which is always a
transitivity set) and its complement. On the other hand, it is easy to see that the last condition is 
equivalent to the doubly transitivity of $G(P).$ 

Furthermore, the Riemann-Hurwitz formula implies that if $h_{P}(x,y)$ is irreducible then
$$
2-2g=\left(
\sum_{i=1}^{r}
\sum_{j_1=1}^{u_{i}} \sum_{j_2=1}^{v_{i}} \GCD(p_{i,j_1}p_{i,j_2})-\mu \right)-(r-2)(n^2-n),
$$
where $\mu$ is the total number of disjointed cycles of permutations  
$\delta_i,$ $1\leq i \leq r,$ on $\Delta$. Since $\mu$ coincides with 
the total number of disjointed cycles of permutations  
$\alpha_i,$ $1\leq i \leq r,$ using the Riemann-Hurwitz formula again 
we see that
$\mu=2+(r-2)n$ and therefore \eqref{rh2} holds. \end{proof}

\section{Irreducibility of $h_{P,Q}(x,y)$ and $h_{P}(x,y)$}
\subsection{Irreducibility of $h_{P,Q}(x,y)$}
\bp \la{prop} Let $P,Q$ be rational functions, $\deg P=n,$ $\deg Q=m.$ 
Then any of the conditions below implies the irreducibility of the curve $h_{P,Q}(x,y)=0$.
\vskip 0.2cm 
\noindent 1) $\f C(P)\cap \f C(Q)$ contains at most one element,   
\vskip 0.2cm 
\noindent 2) $\GCD(n,m)=1$,
\vskip 0.2cm 
\noindent 3) $P$ is a polynomial and $Q$ is a rational function with no multiple poles. 
 
\ep

\begin{proof} Suppose that 1) holds. Without loss of generality we may assume that $\f C(P)\cap \f C(Q)=z_1$ (if $\f C(P)\cap \f C(Q)=\emptyset$ the proof is similar) and that for some $s$, $2\leq s <r,$ the following condition holds: for $i,$  
$2\leq i \leq s,$ the point $z_i$ is a critical value of $P$ but is not a critical value of $Q$ while  
for $i,$ $s< i \leq r,$ the point $z_i$ is a critical value of $Q$ but is not a critical value of $P.$
This implies that for $i,$  
$2\leq i \leq s,$ the permutation $\delta_i$ permutes rows of $M$ in    
accordance with the permutation $\alpha_i$ but transforms each column of $M$ to itself. 
Similarly, for  
$i,$ $s< i \leq r,$ the permutation $\delta_i$ permutes columns of $M$ in    
accordance with the permutation $\beta_i$ but transforms each row of $M$ to itself.

Since by \eqref{per} the permutation $\alpha_1$ is contained in the group generated by 
$\alpha_2,\alpha_3,$ $ ... ,\alpha_r$ the last group is transitive on the set $P^{-1}\{z_0\}$. This implies that the subgroup $\Gamma_1$ of $\Gamma(P,Q)$ generated by $\delta_2,\delta_3, ... ,\delta_s$ acts transitively on the set of rows.
Similarly, the subgroup $\Gamma_2$ of $\Gamma(P,Q)$
generated by $\delta_{s+1},\delta_{s+2}, ... ,\delta_r$ acts transitively on the set of columns. If now $c_{i_1,j_1}$ and $c_{i_2,j_2}$ are two elements
of $M$ and $\gamma_1\in \Gamma_1$ (resp. $\gamma_2\in \Gamma_2$) is an element
such that $i_1^{\gamma_1}=i_2$ (resp. $j_1^{\gamma_2}=j_2$) then 
$$(c_{i_1,j_1})^{\gamma_1\gamma_2}=(c_{i_2,j_1})^{\gamma_2}=c_{i_2,j_2}.$$
Therefore, the subgroup of $\Gamma(P,Q)$
generated by $\delta_{2},\delta_3, ... ,\delta_r$ acts transitively on the set of elements of $M$ and hence 
the action of the group $\Gamma(P,Q)$  is also transitive.

In order to prove the sufficiency of 2) it is enough to observe that since
for any $j,$ $1\leq j \leq o(P,Q),$
the imprimitivity system $\Omega_P(j)$
(resp. $\Omega_Q(j)$) contains $n$ (resp. $m$ blocks), 
the cardinality of any set $U_j,$ $1\leq j \leq o(P,Q),$ is divisible by the $\LCM(n,m)$. 
On the other hand, if 2) holds then $\LCM(n,m)=mn$. Since $M$ contains $mn$ elements this implies that the group $\Gamma(P,Q)$ is transitive. 

Suppose finally 
that 3) holds. Without loss of generality we may assume that $z_1=\infty$. 
Let $c_{i_1,j_2}$ and $c_{i_2,j_2}$ be two elements of $M.$  
Since the group $\beta_1,\beta_1,$ $ ... ,\beta_r$ is transitive on the set $Q^{-1}\{z_0\}$
there exists $g\in \Gamma(P,Q)$ such that 
$(c_{i_1,j_1})^{g}=c_{i,j_2}$ for some $i,$ $1\leq i \leq n$. On the other hand, since $P$ is a polynomial the permutation $\alpha_1$ is 
a full cycle and hence there exists a number $k,$ $1\leq k \leq n,$ such that $i^{\delta_1^k}=i_2$. Furthermore, since $Q$ has no multiple poles
the permutation $\delta_1$ transforms each column of $M$ to itself. Therefore, $$(c_{i_1,j_1})^{g\delta_1^k}=(c_{i,j_1})^{\delta_1^k}=c_{i_2,j_2}$$
and hence the group $\Gamma(P,Q)$ is transitive. \end{proof}

If rational functions $P$ and $Q$ have two common critical values then the curve $h_{P,Q}(x,y)$ can be reducible. Nevertheless, it turns out that all
reducible curves $h_{P,Q}(x,y)$ for which $\f C(P)\cap \f C(Q)$ contains two elements can be 
described explicitly. 
In order to obtain such a description (and another proof of the first part of Proposition \ref{prop}) we will use the following result which is due to Fried (see \cite{f4}, Proposition 2, \cite{f2}, Lemma 4.3, or \cite{pak}, Theorem 3.5).

For a rational function $F=F_1/F_2$ denote by $\Omega_F$ the splitting field 
of the polynomial $F_1(x)-zF_2(x)=0$ over $\C(z).$ 

\bp \la{frr} {\rm (\cite{f4})} Let $P,Q$ be rational functions such that the curve \linebreak $h_{P,Q}(x,y)$ is reducible. Then there exist rational functions 
$A,B,\tilde P, \tilde Q$ such that 
\be \la{lk} P=A \circ \tilde P,\ \ \ \ Q=B \circ \tilde Q, \ \ \ \ o(A,B)=o(P,Q),  \ \ \ \ \Omega_A=\Omega_B.\ee
In particular, it follows from $\Omega_A=\Omega_B$ that
$\f C(A)=\f C(B).$ \qed
\ep

Notice that since for the functions $A,B$ in Proposition \ref{frr} the inequality $o(A,B)=o(P,Q)>1$ holds the degrees of $A,B$ are greater than 1. 

The proposition below supplements the first part of Proposition \ref{prop}.

\bp \la{ebs} Let $P,Q$ be rational functions such that $\f C(P)\cap \f C(Q)$ contains two elements. Then 
the curve $h_{P,Q}(x,y)$ is reducible if and only if there exist rational functions $P_1, Q_1$ and 
a M\"obius transformation $\mu$ 
such that 
\be \la{for} P=\mu\circ z^d \circ P_1, \ \ \ \ Q=\mu\circ z^{d} \circ Q_1 \ee 
for some integer $d>1$. \ep 

\begin{proof}  
Suppose that $h_{P,Q}(x,y)$ is reducible and let $A,B,\tilde P, \tilde Q$ be rational functions 
from Proposition \ref{frr}. Set $C=\f C(A)=\f C(B).$
By the chain rule 
$$\f C(P)=\f C(A)\cup A(\f C(\tilde P)), \ \ \ \ \f C(Q)=\f C(B)\cup B(\f C(\tilde Q))$$ 
and therefore $C\subseteq \f C(P)\cap \f C(Q)$. 
Therefore, since $\card\{\f C(P)\cap \f C(Q)\}=2$ 
and the degrees of $A,B$ are greater than 1, each of the functions $A$ and $B$ has exactly two critical values.

It follows from equality \eqref{per} that 
for the permutations $\kappa_1,\kappa_2$ corresponding to the critical values of $A$
the equality $\kappa_1\kappa_2=1$ holds.  
Therefore each of these permutations is a cycle of length $d=\deg A$ and this implies easily that
there exist M\"obius transformations $\mu$ and $\nu$ such that 
$A=\mu \circ z^{d} \circ \nu.$ Similarly, $B=\tilde \mu \circ z^{\tilde d}\circ \tilde \nu$ for some M\"obius transformations $\tilde \mu, \tilde \nu$ and $\tilde d=\deg B.$
Furthermore, it follows from $\Omega_A=\Omega_B$ that $\tilde d=d$ 
and the equality $\f C(A)=\f C(B)$ implies that $\tilde \mu=\mu \circ cz^{\pm 1}$ 
for some $c\in \C$. Setting now $$P_1=\nu \circ \tilde P,\ \ \ Q_1=c^{1/d}z^{\pm 1}\circ \tilde \nu \circ \tilde Q$$
we conclude that \eqref{for} holds for some $d>1$.

Finally, it is clear that if \eqref{for} holds then the curve $h_{P,Q}(x,y)$ is reducible. \end{proof}

\subsection{Irreducibility of $h_{P}(x,y)$}
Recall that a rational function $P$ is called decomposable if 
there exist rational functions $P_1,P_2,$ $\deg P_1>1,$ $\deg P_2 >1,$ such that 
$P=P_1\circ P_2.$ Otherwise, $P$ is called indecomposable.

It is easy to see that if the curve $h_P(x,y)$ is irreducible then $P$ is necessarily indecomposable. Indeed,  
since the curve $h_{P_1,P_1}(x,y)=0$ has the factor $x-y$, the curve $h_{P_1\circ P_2,P_1\circ P_2}(x,y)=0$
has the factor $h_{P_2,P_2}(x,y)=0$ and hence  
the curve $h_{P_1\circ P_2}(x,y)$ has the factor 
$h_{P_2}(x,y)$.

\bp \la{wi} Let $P$ be an indecomposable rational function. Suppose that $P$ has at least one simple critical value.
Then the curve $h_P(x,y)$ is irreducible.
\ep
\begin{proof}  Indeed, a rational function $P$ is indecomposable if and only if its monodromy group $G(P)$ is primitive. 
Furthermore, if $P$ has a simple critical value $z_j,$ $1\leq j \leq r,$
then the permutation $\alpha_j$ which corresponds to this critical value is a transposition. 
On the other hand, it is known (see e.g. Theorem 13.3 of \cite{wi}) that a primitive permutation group containing 
a transposition is a full symmetric group. Since a symmetric group is doubly transitive Proposition \ref{wi} follows now from Corollary \ref{doub}. \end{proof}

Recall that a point 
$y\in\C\P^1$ is called a critical point of a rational function $P$ if the local multiplicity of $P$ at $y$ is greater than 1. 
Say that a rational function $P$ satisfies the separation condition if for any distinct critical points 
$y_1$, $y_2$ of $P$ the inequality $P(y_1)\neq P(y_2)$ holds. Notice that this condition is often assumed in 
the papers about uniqueness polynomials for meromorphic functions (see e.g. \cite{a0}, \cite{a}, \cite{fu1}, \cite{fu2}, \cite{fu3}). The Proposition \ref{wi2} below shows that the separation condition actually is closely related to the indecomposability condition.

\bp \la{wi2} Let $P$ be a rational function satisfying the separation condition.
Then either $P$ is indecomposable or \be \la{rtf} P=\gamma_1 \circ z^{n} \circ \gamma_2\ee for some M\"obius transformations
$\gamma_1, \gamma_2$ and a composite number $n$. In particular, if $P$ has at least one simple critical value 
then the curve $h_P(x,y)$ is irreducible.
\ep
\begin{proof} First of all observe that for any finite set $T\subset\C\P^1$ and any
rational function $F$ of degree $n$ 
the Riemann-Hurwitz formula implies that
\be \la{ngur} \card\{ F^{-1}\{T\}\}\geq 2+(\card\{T\}-2)n\ee
and the equality attains if and only if $T=\f C(F).$ In particular, 
if $n>1$ then $\card\{ F^{-1}\{T\}\}>\card\{T\}$ unless $T=\f C(F)$ and $\card\{\f C(F)\}=2$. Recall that as it was 
noted in the proof of Proposition \ref{ebs} the equality $\card\{\f C(F)\}=2$ implies that there exist M\"obius transformations $\mu$ and $\nu$ such that 
$F=\mu \circ z^{n} \circ \nu.$

Suppose now that a rational function $P$ satisfying the separation condition is decomposable and let $P_1,$ $P_2$ be rational functions of degree greater than 1 such that $P=P_1\circ P_2$. Denote by $\f S(P_1)$ the set of critical points of $P_1.$ It follows from the chain rule that 
if $\zeta\in \f S(P_1)$ then any 
point $\mu$ such that $P_2(\mu)=\zeta$ is a critical 
point of $P$. 
Therefore, the separation condition implies that for any $\zeta\in \f S(P_1)$ the set $P^{-1}_2\{\zeta\}$ consists of a unique point and hence 
\be \la{pizd} \card\{ P_2^{-1}\{\f S(P_1)\}\} =\card\{\f S(P_1)\} .\ee
As it was observed above \eqref{pizd} implies that $\f S(P_1)=\f C(P_2)$, $\card\{\f C(P_2)\}=2$, and $P_2=\gamma_2 \circ 
z^{d_2} \circ\alpha_2$ for some M\"obius transformations $\alpha_2, \gamma_2$ and $d_2>1.$ 

Furthermore, it follows from $\card\{\f S(P_1)\}=2$ that $\card\{\f C(P_1)\}=2$ and therefore $P_1=\alpha_1 \circ z^{d_1} \circ \gamma_1$ 
for some M\"obius transformations $\alpha_1, \gamma_1$ and $d_1>1.$ 
Since $\f S(P_1)=\f C(P_2)$ we have $\gamma_1\circ \gamma_2=\pm z$ and hence
\eqref{rtf} holds for $\gamma_1=\alpha_1,$ $\gamma_2= z^{\pm 1} \circ \alpha_2$, $n=d_1d_2.$ 
Finally, if $P$ has at least one simple critical value then it may not have the form \eqref{rtf} and hence the curve $h_P(x,y)$ is irreducible by Proposition \ref{wi}.
\end{proof}

\bc \la{wi3} Let $P$ be a rational function which has only simple critical values.
Then the curve $h_P(x,y)$ is irreducible. 
\ec
\begin{proof}  Indeed, a critical value $\zeta$ of a rational function $P$ is simple if and only if the set $P^{-1}\{\zeta\}$ contains a unique critical point and the local multiplicity of $P$ at this point is 2. Therefore, 
if $P$ has only simple critical values then $P$
satisfies the separation condition and hence $h_P(x,y)$ is irreducible by Proposition \ref{wi2}.\end{proof}

\section{Equations $P\circ f=Q\circ g$ and $P\circ f=cP\circ g$}
\subsection{Equation $P\circ f=Q\circ g$}

\noindent{\it Proof of Theorem \ref{1}.} Since $\f C(P)\cap \f C(Q)=\emptyset$ it follows from the first part of Proposition \ref{prop} that the curve $h_{P,Q}(x,y)=0$
is irreducible. Therefore, in view of the Picard theorem in order to prove the theorem it is enough to check that the genus of $h_{P,Q}(x,y)=0$ equals $(n-1)(m-1)$.

We will keep the notation of section 2. Without lost of generality we may assume that there exists $s$, $1< s <r,$ such that for $i,$ $1\leq i \leq s,$ the point $z_i$ 
is a critical value of $P$ but is not a critical value of $Q$ while  
for $i,$ $s< i \leq r,$ the point $z_i$ is a critical value of $Q$ but is not a critical value of $P.$
Then by Corollary \ref{p3} we have:

$$ 2-2g=
\sum_{i=1}^{r}
\sum_{j_1=1}^{u_{i}} \sum_{j_2=1}^{v_{i}} \GCD(p_{i,j_1}q_{i,j_2})-(r-2)nm=$$
$$=\sum_{i=1}^{s}
\sum_{j_1=1}^{u_{i}} \sum_{j_2=1}^{v_{i}} \GCD(p_{i,j_1}q_{i,j_2})+\sum_{i=s+1}^{r}
\sum_{j_2=1}^{v_{i}}\sum_{j_1=1}^{u_{i}}  \GCD(p_{i,j_1}q_{i,j_2})-(r-2)nm=
$$
$$=\sum_{i=1}^{s}
\sum_{j_1=1}^{u_{i}} \sum_{j_2=1}^{m} 1+ \sum_{i=1}^{s}
\sum_{j_2=1}^{v_{i}} \sum_{j_1=1}^{n} 1 -(r-2)nm=
\sum_{i=1}^{s}
\sum_{j_1=1}^{u_{i}} m+ \sum_{i=1}^{s}\sum_{j_2=1}^{v_{i}} n-(r-2)nm.$$
Since by the Riemann-Hurwitz formula we 
have:
$$\sum_{i=1}^{s}
\sum_{j_1=1}^{u_{i}}1=(s-2)n+2, \ \ \  \sum_{i=1}^{s}\sum_{j_2=1}^{v_{i}}1=(r-s-2)m+2,$$
this implies that 
$$2-2g=((s-2)n+2)m+((r-s-2)m+2)n-(r-2)nm=2m+2n-2mn.
$$ 
Therefore, $$g=nm -m-n +1=(m-1)(n-1).\ \ \Box$$ \vskip 0.2cm

\noindent{\it Proof of Theorem \ref{ca1}.} First of all remove from  
$\C\P^{2n+1}\times\C\P^{2m+1}$ the hyperplanes $b_n=0$ and $d_m=0$. Then we may 
set $b_n=1,$ $d_m=1$ and to identify the pair $P,Q$ 
with the point $(a_n,...,a_0,b_{n-1}, ... , b_0,
c_m,...,c_0,d_{m-1}, ... , d_0)$ of the affine space $\C^{2n+2m+2}.$
Notice that the condition $b_n\neq 0,$ $d_m\neq 0$ implies that 
the point $\infty$ can not be a critical point of $P$ or $Q$ corresponding to the critical value $\infty$. 
Furthermore, remove from $\C^{2n+2m+2}$ the hyperplanes $\Gamma_1$ and $\Lambda_1$ corresponding to the discriminants of the polynomials
$$B(z)=z^n+b_{n-1}z^{n-1}+...+b_0, \ \ \ \ D(z)=z^m+d_{m-1}z^{m-1}+...+d_0.$$ 
Then for remaining pairs $P,Q$ the finite points from $\C\P^1$ also can not be critical points 
corresponding to the critical value $\infty$. 
Finally, remove the hyperplanes $\Gamma_2 :\,a_{n-1}-b_{n-1}a_n=0$ and 
$\Lambda_2:\, c_{m-1}-d_{m-1}c_m=0$ containing functions for 
which the point $\infty$ is a critical point.
If now $P,Q$ is a pair from $\C^{2n+2m+2}\setminus \Gamma$, where  $\Gamma=\Gamma_1\cup  \Gamma_2\cup \Lambda_1\cup \Lambda_2$, then all critical values and critical points of $P,Q$ are finite.

Set $$E(z)=A^{\prime}(z)B(z)-A(z)B^{\prime}(z),\ \ \  F(z)=C^{\prime}(z)D(z)-C(z)D^{\prime}(z),$$
where $$A(z)=a_nz^n+a_{n-1}z^{n-1}+...+a_0,\ \ \ C(z)=c_mz^m+c_{m-1}z^{m-1}+...+c_0.$$ 
By construction, if $P,Q$ is a pair from $\C^{2n+2m+2}\setminus \Gamma$ then any critical point of $P$ (resp. of $Q$) is a zero of the polynomial $E$ (resp. of $F$). Furthermore, the set of critical values of $P$ (resp. of $Q$) coincides with the set of zeros of the polynomial $U(x)$ (resp. of the polynomial $V(x)$), where
$$U(x)=\Res_z(E(z),A(z)-xB(z))\ \ \ V(x)=\Res_z(F(z),C(z)-xD(z)),$$ 
and the corresponding resultants are considered as polynomials in $x$.
Therefore, after removing from $\C^{2n+2m+2}\setminus \Gamma$
the hyperplane corresponding to 
$$\Res_x(U(x),V(x))$$ all remaining pairs $P,Q$ have different critical values and corollary follows
from Theorem A.
$\Box$ \vskip 0.2cm

Clearly, using formula \eqref{rh0} one can obtain other criteria, similar to Theorem \ref{1},
for equation \eqref{ur} to have only trivial solutions.
However, the finding of a complete list of rational functions for which the curve $h_{P,Q}(x,y)$ has a factor of genus 0 or 1,
or equivalently the equation $P\circ g=Q\circ g$ has non-constant meromorphic solutions, seems to be a very difficult problem. 
Let us mention several particular cases when the answer is known.

If $P,Q$ are polynomials then the description of curves $h_{P,Q}(x,y)$ having a factor of genus zero
with one point at infinity is equivalent to the 
classification of polynomial solutions of the equation \be \la{ta} P\circ F=Q\circ G.\ee
The last problem was essentially solved by Ritt in his classical paper \cite{ri}.
Notice that equation \eqref{ta} is closely connected with the problem of description of
polynomials $F,G$ satisfying the equality $F^{-1}\{S\}=G^{-1}\{T\}$ for some compact sets $S, T\subset \C$ (see \cite{p}).

A more general question 
of description of curves $h_{P,Q}(x,y)$ having a factor
of genus $0$ with at most two points at infinity is related with the number theory
and was studied in the papers of Fried \cite{f1} and Bilu $\&$ Tichy \cite{bilu}. In particular, in \cite{bilu}
an explicit list of such curves was obtained. 
Finally, the classification of solutions of the equation
$$L=A\circ B=C\circ D,$$ where $L$ is a rational function with at most two poles and $A,B,C,D$ are arbitrary rational functions, was obtained
in the recent papers \cite{pak0}, \cite{pak} (see also \cite{mp}).
Notice that this classification, generalizing the Ritt theorem and the classification of Bilu and Tichy, also permits to describe solutions of the functional equation $$h=P(f)=Q(g),$$ where $P,Q$ are rational functions and $f,g,h$ are entire functions (see \cite{pak2}). 
In its turns it gives an explicit description of strong uniqueness
polynomials for entire functions (\cite{pak2}).

Another important result about curves $h_{P,Q}(x,y)$, obtained by Avanzi and Zannier \cite{az}, is the classification of 
polynomials $P$ such that the curve $P(x)-cP(y)=0$ has a factor of genus zero for some $c\in \C.$ Notice that this result solves ``a half" 
of the problem of description of strong uniqueness polynomials for meromorphic functions. However, an extension
of the classification of \cite{az} which would include also factors of genus 1 does not seem to be an easy problem. 

Finally, notice that in the other paper by Avanzi and Zannier \cite{az0} was obtained the classification of 
curves $h_{P,Q}(x,y)$ of genus 1 under condition that $\GCD(\deg P,\deg Q)=1$. Observe that together with the Ritt theorem this gives
a complete classification of polynomials such that $\GCD(\deg P,\deg Q)=1$ and the equation $P\circ f=Q\circ g$ 
has non-constant meromorphic solutions.

\subsection{Equation $P\circ f=cP\circ g$}
\noindent{\it Proof of Theorem \ref{2}.} We will keep the notation of Section 2.
First of all observe that all critical values of a rational function $P,$ $\deg P=n,$ are simple if and only if for   
the number of critical values $r$ of $P$ the equality \be \la{ner} r=2n-2\ee holds. Indeed, if all critical values of $P$ are simple then 
\be \la{ko} \lambda_i=(2,1,1,\dots,1), \ \ \  u_i=n-1, \ \ \ 1\leq i \leq r, \ee and therefore by
the Riemann-Hurwitz formula we have:
\be \la{hr} 2=\sum_{i=1}^r u_i -(r-2)n=2n-r.\ee
On the other hand, if \eqref{ner} holds then the Riemann-Hurwitz formula implies that \be \la{wer} \sum_{i=1}^{2n-2} u_i =2n^2-4n+2=(n-1)(2n-2).\ee Since for any $i,$ $1\leq i \leq 2n-2,$ the inequality 
$u_i \leq n-1$ holds and the equality attains if and only if 
$\lambda_i=(2,1,1,\dots,1)$, it follows from \eqref{wer} that all critical values of $P$ are simple.

Furthermore, by Corollary \ref{wi3} the curve $h_P(x,y)$ is irreducible. Since 
\eqref{ko} implies that for any $i,$ $1\leq i \leq r,$ 
$$
\sum_{j_1=1}^{u_{i}}\sum_{j_2=1}^{v_{i}} \GCD(p_{i,j_1}p_{i,j_2})=n^2-2n+2
$$
it follows from Corollary \ref{doub} taking into account \eqref{ner} that
$$4-2g=
\sum_{i=1}^{r}
\sum_{j_1=1}^{u_{i}} \sum_{j_2=1}^{v_{i}} \GCD(p_{i,j_1}p_{i,j_2})-(r-2)n^2=r(n^2-2n+2)-(r-2)n^2=$$
$$=(2n-2)(n^2-2n+2)-(2n-4)n^2=-2n^2+8n-4. 
$$
Hence $g=(n-2)^2$ and therefore $g$ is less than 2 if and only if $n< 4.$ \qed

\noindent{\it Proof of Theorem \ref{ca2}.} We will keep the notation used in the proof of Theorem \ref{ca1}. 
First of all remove from $\C\P^{2n+1}$ the hyperplane $b_n=0$ and 
identify a rational function $P$ with the point $(a_n,...,a_0,b_{n-1}, ... , b_0)$
of the affine space $\C^{2n+1}.$
Furthermore, remove from $\C^{2n+1}$ the hyperplanes $\Gamma_1$ and $\Gamma_2.$ 
As above  
if $P\in \C^{2n+1}\setminus \{\Gamma_1\cup \Gamma_2\}$ then any critical point of $P$ is a zero of the polynomial 
$E(z)$
and critical values of $P$ coincide with zeros of the polynomial $U(x).$

Furthermore, after removing from $\C^{2n+1}\setminus \{\Gamma_1\cup \Gamma_2\}$ the hyperplane $\Omega_1$ corresponding to the discriminant of the polynomial $U(x)$ any remaining function $P$
has $$\deg_x U=\deg_z E=2n-2$$ distinct critical values. As it was observed in the proof of Theorem \ref{2} this implies that all critical values of $P$ are simple. In particular, by Theorem \ref{2} the curve $h_P(x,y)$ is irreducible
and of genus $>1.$ 

Consider a polynomial in $y$ defined by the expression $$L(y)=\Res_x(U(x),y^{2n-2}U(x/y)).$$ It is easy to see that 
$\deg L(y)=2n-2$ and that the set of zeros of $L(y)$ coincides with the set $C_P$ consisting of numbers $\alpha\in \C^{\ast}$ such that $\f C(P)\cap \f C(\alpha P)\neq\emptyset.$ Furthermore, it follows easily from the definition 
of the resultant that $y=1$ is a root of multiplicity $2n-2$ of $L(y).$ 
Set 
$$W(y)=\frac{L(y)}{(y-1)^{2n-2}}$$ and define $\Omega_2$ as the hyperplane of $\C^{2n+1}$ corresponding to 
the discriminant of $W(y).$

If $P\in \C^{2n+1}\setminus \Omega$, where $\Omega=\{\Gamma_1\cup \Gamma_2\cup \Omega_1\cup \Omega_2\}$ then the set $C_P$ contains  
$$\deg W(y)=(2n-2)^2-(2n-2)=(2n-2)(2n-3)$$ different elements distinct from $1$. On the other hand, 
if $$\f C(P)=\{z_1,z_2,\dots, z_{2n-2}\}$$ then any element $\alpha\in C_P,$ $\alpha\neq 1,$ should have the form 
$z_i/z_j$ for some distinct $i,j,$ $1\leq i,j \leq 2n-2,$ and therefore $C_P\setminus \{1\}$ contains at most $$2C_{2n-2}^2=(2n-2)(2n-3)$$ elements and the equality attains if and only if for any $\alpha\in C_P,$ $\alpha\neq 1,$ the set $\f C(P)\cap \f C(\alpha P)$ contains exactly one element.  

Hence, if $P\in \C^{2n+1}\setminus \Omega$ then for any $c\in \C,$ $c\neq 1,$ the intersection $\f C(P)\cap \f C(cP)$
contains at most one element and therefore the curve $h_{P,cP}(x,y)$ is irreducible by Proposition \ref{prop}. If 
$\f C(P)\cap \f C(cP)=\emptyset$
then by Theorem \ref{1} the genus of $h_{P,cP}(x,y)$ equals $(n-1)^2$. 
On the other hand, if $\f C(P)\cap \f C(cP)$ 
contains a single element then it is easy  
to calculate using formula \eqref{rh0} and taking into account equalities \eqref{ko} that the genus of $h_{P,cP}(x,y)$ equals $n^2-2n.$ 
In both cases the assumption $n\geq 4$ implies that the genus of $h_{P,cP}(x,y)$
is greater than 1.
\qed

\end{document}